\documentclass[draftcls,a4paper,11pt,onecolumn]{IEEEtran}%
\usepackage{amsfonts}
\usepackage{amsmath}
\usepackage[latin1]{inputenc}
\usepackage[normalem]{ulem}
\usepackage{color}
\usepackage{amssymb}
\usepackage{graphicx}%
\usepackage{citesort}
\newcounter{algo}
\renewcommand{\thealgo}{\arabic{algo}}

\graphicspath{{./figures}}
\begin{document}

\title{Bayesian Inference for Linear Dynamic Models with Dirichlet Process Mixtures}
\author{F.\ Caron$^{1}$, M.\ Davy$^{1}$, A.\ Doucet$^{2}$, E.\ Duflos$^{1}$,
P.\ Vanheeghe$^{1}$.\\$^{1}$ CNRS/Ecole Centrale de Lille and INRIA-FUTURS SequeL  team, Villeneuve d'Ascq, France\\$^{2}$ University of British Columbia, Canada}
\date{}
\maketitle

\begin{abstract}
Using Kalman techniques, it is possible to perform optimal estimation in
linear Gaussian state-space models. We address here the case where the noise
probability density functions are of unknown functional form. A flexible
Bayesian nonparametric noise model based on Dirichlet process mixtures is
introduced. Efficient Markov chain Monte Carlo and Sequential Monte Carlo
methods are then developed to perform optimal batch and sequential estimation
in such contexts. The algorithms are applied to blind deconvolution and change
point detection. Experimental results on synthetic and real data demonstrate
the efficiency of this approach in various contexts.

\end{abstract}




\section{Introduction}

Dynamic linear models are used in a variety of applications, ranging from
target tracking, system identification, abrupt change detection, etc. The
models are defined as follows~:%
\begin{eqnarray}
\mathbf{x}_{t}&=&A_{t}\mathbf{x}_{t-1}+C_{t}\mathbf{u}_{t}+G_{t}\mathbf{v}_{t}
\label{eq:statemodel}\\
\mathbf{z}_{t}&=&H_{t}\mathbf{x}_{t}+\mathbf{w}_{t} \label{eq:observationmodel}%
\end{eqnarray}
where $\mathbf{x}_0\sim \mathcal{N}(\mu_0,\Sigma_0)$, $\mathbf{x}_{t}$ is the hidden state vector, $\mathbf{z}_{t}$ is the
observation, $\mathbf{v}_{t}$ and $\mathbf{w}_{t}$ are sequences of mutually
independent random variables such that $\mathbf{v}_t  \stackrel{\text{i.i.d.}}{\sim} F^v$ and $\mathbf{w}_t  \stackrel{\text{i.i.d.}}{\sim} F^w$. $A_{t}$ and $H_{t}$ are the known state and
observation matrices, $\mathbf{u}_{t}$ is a known input, $C_{t}$ the input
transfer matrix and $G_{t}$ is the state transfer matrix. Let us denote
$\mathbf{a}_{i:j}=\left(  \mathbf{a}_{i},\mathbf{a}_{i+1},...,\mathbf{a}%
_{j}\right)  $ for any sequence $\left\{  \mathbf{a}_{t}\right\}  $. The main
use of model (\ref{eq:statemodel})-(\ref{eq:observationmodel}) is to estimate
the hidden state $\mathbf{x}_{t}$ given the observations $\mathbf{z}_{1:t}$
(filtering, with a forward recursion) or $\mathbf{z}_{1:T}$ for $t\leq T$
(smoothing, with a forward-backward recursion).

It is a very common choice to assume that the noise probability density functions (pdfs) $F_v$ and $F_w$ are Gaussian, with known parameters, as this enables the
use of Kalman filtering/smoothing. In such a framework, Kalman techniques are
optimal in the sense of minimizing the mean squared error. There are, however,
a number of cases where the Gaussian assumption is inadequate, e.g. the actual
observation noise distribution or the transition noise are multimodal (in
Section~\ref{sec:applications}, we provide several such examples). In this
paper, we address the problem of \emph{optimal state estimation when the
probability density functions of the noise sequences are unknown and
need to be estimated on-line or off-line from the data}. This problem takes
place in the class of identification/estimation of linear models
with unknown statistic noises.

\subsection{Proposed approach}

\label{sec:proposed-approach} Our methodology\footnote{Preliminary results were presented in Caron et al.~\cite{Caron2006}.} relies on the introduction of a
Dirichlet Process Mixture (DPM), which is used to model the unknown pdfs of
the state noise $\mathbf{v}_{t}$ and measurement noise $\mathbf{w}_{t}$. DPMs
are flexible Bayesian nonparametric models which have become very popular in
statistics over the last few years, to perform nonparametric density
estimation~\cite{Walker99,Neal00,Muller04}. Briefly, a realization of a DPM
can be seen as an \emph{infinite} mixture of pdfs with given parametric shape
(e.g., Gaussian) where each pdf is denoted $f(\cdot|\theta)$. The parameters
of the mixture (mixture weights and locations of the $\theta$'s) are given by
the random mixture distribution $\mathbb{G}(\theta)$,\ which is sampled from a
so-called \emph{Dirichlet Process}. A prior distribution, denoted
${\mathbb{G}}_{0}(\theta)$ must be selected over the $\theta$'s (e.g.,
Normal-Inverse Wishart for the DPM of Gaussians case, where $\theta$ contains
the mean vector and the covariance matrix), while the weights follow a
distribution characterized by a positive real-valued parameter $\alpha$. For
small $\alpha$, only a small fraction of the weights is significantly nonzero,
whereas for large $\alpha$, many weights are away from zero. Thus, the
parameter $\alpha$ tunes the prior distribution of components in the mixture, without
setting a precise number of components. Apart from this implicit, powerful
clustering property, DPMs are computationally very attractive due to the
so-called \emph{Polya urn representation} which enables straightforward
computation of the full conditional distributions associated to the latent
variables $\theta$.

\subsection{Previous works}

\label{sec:previous-works} Several algorithms have been developed to estimate
noise statistics in linear dynamic
systems~\cite{Mehra1970,Myers1976,Maine1981,Maryak2004}. However, these
algorithms assume Gaussian noise pdfs (with unknown mean and covariance
matrix). As will be made clearer in the following, this is a special case of
our framework: if the scaling coefficient $\alpha$ tends to $0$, the
realizations of the DPM of Gaussian pdfs converge in distribution to a single
Gaussian with parameter prior distribution given by the base distribution
$\mathbb{G}_{0}$. Algorithms have also been developed to deal with
non-Gaussian noises distributions, such as student-t~\cite{Shephard1994},
$\alpha$-stable~\cite{Lombardi2005} or mixture of Gaussians~\cite{Carter96}.
These works are based on a given prior parametric shape of the pdf which we do
not assume in this paper.

Though many recent works have been devoted to DPMs in various contexts such as
econometrics~\cite{Griffin04b}, geoscience~\cite{Pievatolo00}\ and
biology~\cite{Do05,Medvedovic02}, this powerful class of models has never been
used in the context of linear dynamic models (to the best of our knowledge).
In this paper, we show that DPM-based dynamic models with unknown noise
distributions can be defined easily. Moreover, we provide several efficient
computational methods to perform Bayesian inference, ranging from Gibbs
sampling (for offline estimation) to Rao-Blackwellized particle filtering for
online estimation.

\subsection{Paper organization}

\label{sec:paper-organization} This paper is organized as follows. In
Section~\ref{sec:DPM}, we recall the basics of Bayesian nonparametric density
estimation with DPMs. In Section~\ref{sec:DynaModelWithUND} we present the
dynamic model with unknown noise distributions. In Section~\ref{sec:MCMC} we
derive an efficient Markov chain Monte Carlo (MCMC) algorithm to perform
optimal estimation in the batch (offline) case. In Section~\ref{sec:SMC}, we
develop a Sequential Monte Carlo (SMC) algorithm/Particle filter to perform
optimal estimation in the sequential (online) case. All these algorithms can
be interpreted as Rao-Blackwellized methods. In Section~\ref{sec:discussion},
we discuss some features of these algorithms, and we relate them to other
existing approaches. Finally, in Section~\ref{sec:applications}, we
demonstrate our algorithms on two applications:\ blind deconvolution of
impulse processes and a change point problem in biomedical
time series. The last section is devoted to conclusions and future research directions.

\section{Bayesian nonparametric density estimation}

\label{sec:DPM}In this section, we review briefly Bayesian nonparametric
density estimation\footnote{There are many ways to understand 'nonparametric'.
In this paper, we follow many other papers in the same
vein~\cite{Walker99,Neal00,Muller04}, where 'nonparametric' refers to the fact
that the pdf of interest cannot be defined by a functional expansion with a
finite-dimensional parameter space.}. We introduce Dirichlet processes as
probabilistic measures on the space of probability measures, and we outline
its discreteness. Then, the DPM model in presented.

\subsection{Density estimation}

Let $\mathbf{y}_{1},...,\mathbf{y}_{n}$ be a statistically exchangeable
sequence distributed with
\begin{equation}
\mathbf{y}_{k}\sim F(\cdot)
\end{equation}
where $\sim$ means \emph{distributed according to}. We are interested here in
estimating $F(\cdot)$ and we consider the following nonparametric model
\begin{equation}
F(\mathbf{y})=\int_{\Theta}f(\mathbf{y}|\theta)d\mathbb{G}(\theta)
\label{eq:integral}%
\end{equation}
where $\theta\in\Theta$ is called the latent variable or cluster variable,
$f(\cdot|\theta)$ is the mixed pdf and $\mathbb{G}(\cdot)$ is the mixing
distribution. Within the Bayesian framework, it is assumed that $\mathbb{G}%
(\cdot)$ is a Random Probability Measure (RPM) \cite{Muller04} distributed
according to a prior distribution (i.e., a distribution over the set of
probability distributions). We will select here the RPM\ to follow a Dirichlet
Process (DP) prior.

\subsection{Dirichlet Processes}

Ferguson \cite{Ferguson73} introduced the Dirichlet Process (DP) as a
probability measure on the space of probability measures. Given a probability
measure $\mathbb{G}_{0}(\cdot)$ on a (measurable)\ space $(\mathcal{T}%
,\mathcal{A})$ and a positive real number $\alpha$, a probability distribution
$\mathbb{G}(\cdot)$ distributed according to a DP of base distribution
$\mathbb{G}_{0}(\cdot)$ and scale factor $\alpha$, denoted $\mathbb{G}%
(\cdot)\sim DP(\mathbb{G}_{0}(\cdot),\alpha)$, satisfies for any partition
$A_{1},...,A_{k}$ of $\mathcal{T}$ and any $k$%
\begin{equation}
\left(  \mathbb{G}(A_{1}),...,\mathbb{G}(A_{k})\right)  \sim\mathcal{D}\left(
\alpha\mathbb{G}_{0}(A_{1}),...,\alpha\mathbb{G}_{0}(A_{k})\right)
\label{eq:defDP}%
\end{equation}
where $\mathcal{D}$ is a standard Dirichlet distribution, classically defined
for a set of random variables $(b_{0},..,b_{p})\sim\mathcal{D(}a_{0}%
,..,a_{p})$ by%
\begin{equation}
\mathcal{D(}a_{0},..,a_{p})=\frac{\Gamma(\sum_{l=0}^{p}a_{l})}{\prod_{l=0}%
^{p}\Gamma(a_{l})}\prod_{l=0}^{p}b_{l}^{a_{l}-1}\delta_{1}(\sum_{l=0}^{p}%
b_{l}) \label{eq:defDirichlet}%
\end{equation}
where $\Gamma$ is the gamma function, and $\delta_{u}(v)$ is the Dirac delta
function, which is zero whenever $v\neq u$. From the definition in
Eq.~(\ref{eq:defDP}), it is easy to show that for every $B\in\mathcal{T}$%
\begin{align}
\mathbb{E}\left[  \mathbb{G}(B)\right]   &  =\mathbb{G}_{0}%
(B)\label{eq:mean_DP}\\
\text{var}\left[  \mathbb{G}(B)\right]   &  =\frac{\mathbb{G}_{0}%
(B)\big(1-\mathbb{G}_{0}(B)\big)}{1+\alpha} \label{eq:var_DP}%
\end{align}

An important property is that the realizations of a Dirichlet process are
\emph{discrete}, with probability one. One can show that $\mathbb{G}$ admits the so-called
\emph{stick-breaking} representation, established by Sethuraman
\cite{Sethuraman94}:%
\begin{equation}
\mathbb{G}(\cdot)\;=\;\sum_{j=1}^{\infty}\pi_{j}\delta_{U_{j}}(\cdot)
\label{eq:stickbreak}%
\end{equation}
with $U_{j}\sim\mathbb{G}_{0}(\cdot)$, $\pi_{j}=\beta_{j}\prod_{l=1}%
^{j-1}(1-\beta_{l})$ and\ $\beta_{j}\sim\mathcal{B}(1,\alpha)$ where
$\mathcal{B}$ denotes the beta distribution.
In the following, we omit $(\cdot)$ in $\mathbb{G}(\cdot)$ and other
distributions, to simplify notations. Using Eq.~(\ref{eq:integral}), it
comes that the following flexible prior model is adopted for the unknown
distribution $F$
\begin{equation}
F(\mathbf{y})=\sum_{j=1}^{\infty}\pi_{j}f(\mathbf{y}|U_{j}). \label{eq:1}%
\end{equation}

Apart from its flexibility, a fundamental motivation to use the DP model is
the simplicity of the posterior update. Let $\theta_{1},\ldots,\theta_{n}$ be
$n$ random samples from $\mathbb{G}$%
\begin{equation}
\theta_{k}|\mathbb{G}\overset{\text{i.i.d.}}{\sim}\mathbb{G}%
\end{equation}
where $\mathbb{G}\sim DP(\mathbb{G}_{0},\alpha)$ then the posterior
distribution of $\mathbb{G}|\theta_{1:n}$ is also a DP%
\begin{equation}
\mathbb{G}|\theta_{1:n}\sim DP(\frac{\alpha}{\alpha+n}\mathbb{G}_{0}+\frac
{1}{\alpha+n}\sum_{k=1}^{n}\delta_{\theta_{k}}\,,\,\alpha+n)
\label{eq:posteriorDP}%
\end{equation}
Moreover, it can be shown that the predictive distribution, computed by
integrating out the RPM $\mathbb{G}$, admits the following \emph{Polya urn}
representation~\cite{Blackwell73}%
\begin{equation}
\theta_{n+1}|\theta_{1:n}\sim\frac{1}{\alpha+n}\sum_{k=1}^{n}\delta
_{\theta_{k}}+\frac{\alpha}{\alpha+n}\mathbb{G}_{0}. \label{eq:Polya}%
\end{equation}
Therefore, conditionally on the latent variables $\theta_{1:n}$ sampled
previously, the probability that a new sample is identical to an existing one
is overall $\frac{n}{\alpha+n}$, whereas, with probability $\frac{\alpha
}{\alpha+n}$, the new sample is distributed (independently) according to
$\mathbb{G}_{0}$. It should be noted that several $\theta_{k}$'s might have
the same value, thus the number of \textquotedblleft alive\textquotedblright%
\ clusters (denoted $M$), that is, the number of distinct values of
$\theta_{k}$, is less than $n$.

The \emph{scaling coefficient} $\alpha$ tunes the number of \textquotedblleft
alive\textquotedblright\ clusters $M$. For large $n$,
Antoniak~\cite{Antoniak74} showed that $\mathbb{E}\left[  M|\alpha,n\right]
\simeq\alpha\log(1+\frac{n}{\alpha})$. As $\alpha$ tends to zero, most of the
samples $\theta_{k}$ share the same value, whereas when $\alpha$ tends to
infinity, the $\theta_{k}$ are almost i.i.d. samples from $\mathbb{G}_{0}$.

\subsection{Dirichlet Process Mixtures}

\label{sec:dirichl-proc-mixt} Using these modeling tools, it is now possible
to reformulate the density estimation problem using the following hierarchical
model known as DPM \cite{Antoniak74}:%
\begin{equation}%
\begin{tabular}
[c]{l}%
$\mathbb{G}\sim DP(\mathbb{G}_{0},\alpha,),$ ~~ and, for $k=1,\ldots,n$\\
~~~~~~~~$\theta_{k}|\mathbb{G}\sim\mathbb{G}$,\\
~~~~~~~~$\mathbf{y}_{k}|\theta_{k}\sim f(\cdot|\theta_{k})$%
\end{tabular}
\ \ \ \ \ \label{eq:DPM}%
\end{equation}
It should be noted that DPMs can model a wide variety of pdfs. In particular,
assuming Gaussian $f(\cdot|\theta_{k})$, the parameter contains both the mean
and the covariance, and, depending on ${\mathbb{G}}_{0}$, the corresponding
DPM may have components with large/small variances.

\subsection{Estimation objectives}

The objective of DPM-based density estimation boils down to estimating the
posterior distribution $p(\theta_{1:n}|\mathbf{y}_{1:n})$, because the
probability $\mathbb{G}$ can be integrated out analytically by using the Polya
urn representation. Although DPMs\ were introduced in the 70's, these models
were too complex to handle numerically before the introduction of Monte Carlo
simulation based methods. Efficient MCMC
algorithms~\cite{Escobar92,Escobar95,Walker99,Neal00,MacEachern00} as well as
Sequential Importance Sampling~\cite{MacEachern99,Fearnhead04} enable to
sample from $p(\theta_{1:n}|\mathbf{y}_{1:n})$. However, these algorithms
cannot be applied to our class of models, which is presented below, because
the noise sequences $\mathbf{v}_{t}$ and $\mathbf{w}_{t}$ are not observed directly.

\section{Dynamic Linear Model with Unknown Noise Distribution}

\label{sec:DynaModelWithUND} The linear dynamic model defined in Eq.'s
(\ref{eq:statemodel})-(\ref{eq:observationmodel}) relies on the unknown noises
$\left\{  \mathbf{v}_{t}\right\}  $ and $\left\{  \mathbf{w}_{t}\right\}  $
distributions, which are assumed to be DPMs in this paper.

\subsection{DPM noise models}

\label{sec:dpm-noise-models} For both $\left\{  \mathbf{v}_{t}\right\}  $ and
$\left\{  \mathbf{w}_{t}\right\}  $, the pdf $f(\cdot|\theta)$ is assumed here
to be a Gaussian, denoted $\mathcal{N}(\mu_{t}^{v},\Sigma_{t}^{v})$ and
$\mathcal{N}(\mu_{t}^{w},\Sigma_{t}^{w})$ respectively. The base distributions
$\mathbb{G}_{0}^{v}$ and $\mathbb{G}_{0}^{w}$ are assumed to be normal inverse
Wishart distributions \cite{Gelman95} denoted $\mathbb{G}_{0}^{v}%
=\mathcal{NI}W(\mu_{0}^{v},\kappa_{0}^{v},\nu_{0}^{v},\Lambda_{0}^{v})$ and
$\mathbb{G}_{0}^{w}=\mathcal{NI}W(\mu_{0}^{w},\kappa_{0}^{w},\nu_{0}%
^{w},\Lambda_{0}^{w})$. The hyperparameters $\psi^{v}$ $=$ $\{\mu_{0}^{v}$,
$\kappa_{0}^{v}$, $\nu_{0}^{v}$, $\Lambda_{0}^{v}\}$ and $\psi^{w}=\{\mu
_{0}^{w},\kappa_{0}^{w},\nu_{0}^{w},\Lambda_{0}^{w}\}$ are assumed fixed but
unknown. Finally, the scale parameters $\alpha^{v}$ and $\alpha^{w}$ are also
assumed fixed and unknown. Overall, the sets of hyperparameters are denoted
$\phi^{v}=\{\alpha^{v},\psi^{v}\}$, $\phi^{w}=\{\alpha^{w},\psi^{w}\}$ and
$\phi=\{\phi^{v},\phi^{w}\}$. For the sake of presentation clarity, we assume
that these hyperparameters are known, but in Subsection~\ref{sec:unkn-hyperp},
we address the case of unknown hyperparameters by defining priors and a
specific estimation procedure.

To summarize, we have the following models%
\begin{equation}%
\begin{tabular}
[c]{l}%
$\mathbb{G}^{v}|\phi^{v}\sim DP(\mathbb{G}_{0}^{v},\alpha^{v}),$
\end{tabular}
\ \ \ \ \ \ \
\begin{tabular}
[c]{l}%
$\mathbb{G}^{w}|\phi^{w}\sim DP(\mathbb{G}_{0}^{w},\alpha^{w}),$
\end{tabular}
\ \ \ \
\end{equation}
and for $t=1,2,\ldots$
\begin{equation}%
\begin{tabular}
[c]{l}%
$\theta_{t}^{v}|\mathbb{G}^{v}%
\overset{\text{i.i.d.}}{\sim}\mathbb{G}^{v},$\\
$\mathbf{v}_{t}|\theta_{t}^{v}\overset{\text{i.i.d.}}{\sim}\mathcal{N}(\mu
_{t}^{v},\Sigma_{t}^{v}).$%
\end{tabular}
\ \ \ \ \ \ \
\begin{tabular}
[c]{l}%
$\theta_{t}^{w}|\mathbb{G}^{w}\overset{\text{i.i.d.}}{\sim}\mathbb{G}^{w},$\\
$\mathbf{w}_{t}|\theta_{t}^{w}\overset{\text{i.i.d.}}{\sim}\mathcal{N}(\mu
_{t}^{w},\Sigma_{t}^{w}).$%
\end{tabular}
\ \ \ \
\end{equation}
where $\theta_{t}^{v}=\{\mu_{t}^{v},\Sigma_{t}^{v}\}\ ($resp. $\theta_{t}%
^{w}=\{\mu_{t}^{w},\Sigma_{t}^{w}\}$) is the latent cluster variable giving
the mean and covariance matrix for that cluster, and $\theta_{t}=\{\theta
_{t}^{v},\theta_{t}^{w}\}$. This model is written equivalently as
$\mathbf{v}_{t}\sim F^{v}(\mathbf{v}_{t})$ and $\mathbf{w}_{t}\sim
F^{w}(\mathbf{w}_{t})$ where $F^{v}$ and $F^{w}$ are fixed but unknown
distributions written as%
\begin{align}
F^{v}(\mathbf{v}_{t})  &  =\int\mathcal{N}(\mathbf{v}_{t};\mu,\Sigma
)d\mathbb{G}^{v}(\mu,\Sigma),\label{eq:defGaussianDPM}\\
F^{w}(\mathbf{w}_{t})  &  =\int\mathcal{N}(\mathbf{w}_{t};\mu,\Sigma
)d\mathbb{G}^{w}(\mu,\Sigma)
\end{align}
In other words, $F^{v}$ and $F^{w}$ are countable infinite mixtures of
Gaussian pdfs of unknown parameters, and the mixing distributions
$\mathbb{G}^{v}$ and $\mathbb{G}^{w}$ are sampled from Dirichlet processes.

\subsection{Estimation of the state parameters}

\label{sec:mixt-param-estim} In this work, our objective is to estimate
${\mathbb{G}}^{v}$ and ${\mathbb{G}}^{w}$ as well as the latent variables
$\{\theta_{t}\}$ and state variable $\{{\mathbf{x}}_{t}\}$ at each time $t$,
conditional on the observations $\{{\mathbf{z}}_{t}\}$. In practice, only the
state variable is of interest -- ${\mathbb{G}}^{v}$, ${\mathbb{G}}^{w}$ and
$\{\theta_{t}\}$ are \emph{nuisance} parameters. Ideally, one would like to
estimate online the sequence of posterior distributions $p(\mathbf{x}%
_{0:t}|\mathbf{z}_{1:t},\phi)$ as $t$ increases or the offline posterior
$p(\mathbf{x}_{0:T}|\mathbf{z}_{1:T},\phi)$, where $T$ is the fixed length of
the observation sequence $\mathbf{z}_{1:T}$. Thanks to the Polya urn
representation, it is possible to integrate out analytically ${\mathbb{G}}%
^{v}$ and ${\mathbb{G}}^{w}$ from these posteriors. The parameters
$\theta_{1:t}$ and $\theta_{1:T}$ remain and the inference is based upon
$p(\mathbf{x}_{0:t},\theta_{1:t}|\mathbf{z}_{1:t},\phi)$ or $p(\mathbf{x}%
_{0:T},\theta_{1:T}|\mathbf{z}_{1:T},\phi)$. The posterior $p(\mathbf{x}%
_{0:t},\theta_{1:t}|\mathbf{z}_{1:t},\phi)$ satisfies for any $t$%
\begin{equation}
p(\mathbf{x}_{0:t},\theta_{1:t}|\mathbf{z}_{1:t},\phi)=p(\mathbf{x}%
_{0:t}|\theta_{1:t},\mathbf{z}_{1:t},\phi)p(\theta_{1:t}|\mathbf{z}_{1:t}%
,\phi). \label{eq:RaoB2}%
\end{equation}

Conditional upon $\theta_{t}$, Eq.'s (\ref{eq:statemodel}%
)-(\ref{eq:observationmodel})\ may be rewritten as
\begin{align}
\mathbf{x}_{t}  &  =F_{t}\mathbf{x}_{t-1}+\mathbf{u}_{t}^{\prime}(\theta
_{t})+G_{t}\mathbf{v}_{t}^{\prime}(\theta_{t})\\
\mathbf{z}_{t}  &  =H_{t}\mathbf{x}_{t}+\mu_{t}^{w}+\mathbf{w}_{t}^{\prime
}(\theta_{t})
\end{align}
where $\mathbf{u}_{t}^{\prime}(\theta_{t})=C_{t}\mathbf{u}_{t}+G_{t}\mu
_{t}^{v}$ and $\mu_{t}^{w}$ are known inputs, $\mathbf{v}_{t}^{\prime}%
(\theta_{t})$ and $\mathbf{w}_{t}^{\prime}(\theta_{t})$ are centered white
Gaussian noise of known covariance matrices $\Sigma_{t}^{v}$ and $\Sigma
_{t}^{w}$, respectively. Thus $p(\mathbf{x}_{0:t}|\theta_{1:t},\mathbf{z}%
_{1:t},\phi)$ (resp. $p(\mathbf{x}_{0:T}|\theta_{1:T},\mathbf{z}_{1:T},\phi)$)
is a Gaussian distribution whose parameters can be computed using a Kalman
filter (resp. smoother)~\cite{Anderson79} for given $\theta_{1:t}$
(resp.$\theta_{1:T}$).

One is generally interested in computing the marginal MMSE state estimate
$\widehat{\mathbf{x}}_{t|t'}^{\text{{\tiny MMSE}}}=\mathbb{E}\left[
\mathbf{x}_{t}|\mathbf{z}_{1:t'}\right] $ (with $t'=t$ or $t'=T$)
\begin{equation}%
\begin{array}
[c]{lcl}%
\widehat{\mathbf{x}}_{t|t'}^{\text{{\tiny MMSE}}} & = & \displaystyle\int
{\mathbf{x}}_{t}p(\mathbf{x}_{t},\theta_{1:t'}|\mathbf{z}_{1:t'},\phi
)d({\mathbf{x}}_{t},\theta_{1:t'})\\
& = & \displaystyle\int{\mathbf{x}}_{t}p(\mathbf{x}_{t}|\theta_{1:t'}%
,\mathbf{z}_{1:t'},\phi)p(\theta_{1:t'}|\mathbf{z}_{1:t'},\phi)d({\mathbf{x}}%
_{t},\theta_{1:t'})\\
& = & \displaystyle\int\widehat{\mathbf{x}}_{t|t'}(\theta_{1:t'})p(\theta
_{1:t'}|\mathbf{z}_{1:t'},\phi)d\theta_{1:t'}\label{eq:3}%
\end{array}
\end{equation}%
where $\widehat{\mathbf{x}}_{t|t}(\theta_{1:t})$ (resp. $\widehat{\mathbf{x}%
}_{t|T}(\theta_{1:T})$) is the mean of the Gaussian $p(\mathbf{x}_{t}%
|\theta_{1:t},\mathbf{z}_{1:t},\phi)$ (resp. $p(\mathbf{x}_{t}|\theta
_{1:T},\mathbf{z}_{1:T},\phi)$). Both $\widehat{\mathbf{x}}_{t|t}(\theta
_{1:t})$ and $\widehat{\mathbf{x}}_{t|T}(\theta_{1:t})$ are computed by the
Kalman filter/smoother, see Sections~\ref{sec:MCMC} and \ref{sec:SMC} below.

Computing these estimates still requires integration w.r.t. the $\theta$'s,
see Eq.~(\ref{eq:3}). This kind of integral is not feasible
in closed-form, but it can be computed numerically by using Monte Carlo
integration~\cite{Robert99}. Briefly, assume that a set of $N$ weighted
samples $\{{\theta}_{1:t}^{(i)}\}_{i=1,\ldots,N}$ with weights ${w}_{t}^{(i)}$
are distributed according to $p(\theta_{1:t}|\mathbf{z}_{1:t},\phi)$, then
e.g., $\widehat{\mathbf{x}}^{\text{{\tiny MMSE}}}_{t|t}$ is computed as
\begin{equation}
\widehat{\mathbf{x}}^{\text{{\tiny MMSE}}}_{t|t}\;\approx\;\sum_{i=1}^{N}%
{w}_{t}^{(i)} \widehat{\mathbf{x}}_{t|t}\big({\theta}_{1:t}^{(i)}%
\big) \label{eq:2}%
\end{equation}
In Eq.~\eqref{eq:2}, the main difficulty consists of generating the weighted
samples $\{{\theta}_{1:t}^{(i)}\}_{i=1,\ldots,N}$ from the marginal posterior
$p(\theta_{1:t}|\mathbf{z}_{1:t},\phi)$ (and similarly, from $p(\theta
_{1:T}|\mathbf{z}_{1:T},\phi)$ in the offline case).

\begin{itemize}
\item For \textbf{offline (batch) estimation $(t=T)$}, this can be done by
MCMC by building a Markov chain of samples $\{{\theta}_{1:T}^{(i)}%
\}_{i=1,\ldots,N}$ with target distribution $p(\theta_{1:T}|\mathbf{z}%
_{1:T},\phi)$ (in that case, ${w}_{t}^{(i)}=1/N$). The
MCMC\ algorithms available in the literature to estimate these Bayesian
nonparametric models -- e.g. \cite{Escobar95,Neal00} -- are devoted to density
estimation in cases where the data are observed directly. They do not apply to
our case because here, the sequences $\left\{  \mathbf{v}_{t}\right\}  $ and
\{$\mathbf{w}_{t}$\} are not observed directly. One only observes $\left\{
\mathbf{z}_{t}\right\}  $, assumed to be generated by the dynamic model
(\ref{eq:statemodel})-(\ref{eq:observationmodel}). Section~\ref{sec:MCMC}
proposes an MCMC algorithm dedicated to this model.

\item For \textbf{online (sequential) estimation}, samples can be generated by
sequential importance sampling, as detailed in Section~\ref{sec:SMC}.
\end{itemize}

\section{MCMC algorithm for off-line state estimation}

\label{sec:MCMC} In this Section, we consider the offline state estimation. As
outlined above, this requires to compute estimates from the posterior
$p(\mathbf{x}_{0:T},\theta_{1:T}|\mathbf{z}_{1:T})$, where we recall that
$\theta_{t}=\{\theta_{t}^{v},\theta_{t}^{w}\}=\{\mu_{t}^{v},\Sigma_{t}^{v}%
,\mu_{t}^{w},\Sigma_{t}^{w}\}$ is the latent variable as defined above. We
first assume that the hyperparameters are fixed and known
(Subsection~\ref{sec:fixed-known-hyperp}), then we let them be unknown, with
given prior distributions (Subsection~\ref{sec:unkn-hyperp}).

\subsection{Fixed and known hyperparameters}

\label{sec:fixed-known-hyperp} In this subsection, the hyperparameter vector
$\phi$ is assumed fixed and known. The marginal posterior $p(\theta
_{1:T}|\mathbf{z}_{1:T},\phi)$ can be approximated through MCMC\ using the
Gibbs sampler~\cite{Robert99} presented in Algorithm 1 below.

\bigskip

\hrule

\smallskip

\noindent\textbf{Algorithm 1:} Gibbs sampler to sample from $p(\theta
_{1:T}|\mathbf{z}_{1:T},\phi)$ \smallskip\hrule\smallskip

\begin{itemize}
\item \underline{Initialization}: For $t=1,...,T$, sample $\theta_{t}^{(1)}$
from an arbitrary initial distribution, e.g. the prior.

\item \underline{Iteration} $i,$ $i=2,\ldots,N^{\prime}+N$:

\begin{itemize}
\item For $t=1,\ldots,T$, sample $\theta_{t}^{(i)}\sim p(\theta_{t}%
|\mathbf{z}_{1:T},\theta_{-t}^{(i)},\phi)$ where $\theta_{-t}^{(i)}%
=\{\theta_{1}^{(i)},..,\theta_{t-1}^{(i)},\theta_{t+1}^{(i-1)},..,\theta
_{T}^{(i-1)}\}$
\end{itemize}
\end{itemize}

\smallskip\hrule$\bigskip$

To implement Algorithm 1, one needs to sample from the conditional pdf
$p(\theta_{t}|\mathbf{z}_{1:T},\theta_{-t},\phi)$ for each of the $N^{\prime
}+N$ iterations (including $N^{\prime}$ burn-in iterations). From Bayes' rule,
we have
\begin{equation}
p(\theta_{t}|\mathbf{z}_{1:T},\theta_{-t},\phi)\propto p(\mathbf{z}%
_{1:T}|\theta_{1:T})p(\theta_{t}|\theta_{-t},\phi). \label{eq:4}%
\end{equation}
where $p(\theta_{t}|\theta_{-t},\phi)=p(\theta_{t}^{v}|\theta_{-t}^{v}%
,\phi^{v})p(\theta_{t}^{w}|\theta_{-t}^{w},\phi^{w})$. From the Polya urn
representation, these two terms are written as (for $w$, replace $v$ with $w$
below):%
\begin{equation}%
\begin{array}
[c]{lcl}%
p(\theta_{t}^{v}|\theta_{-t}^{v},\phi^{v}) & = & \displaystyle \frac{1}%
{\alpha^{v}+T-1}\sum_{k=1,k\neq t}^{T}\delta_{\theta_{k}^{v}}(\theta_{t}%
^{v})\displaystyle +\frac{\alpha^{v}}{\alpha^{v}+T-1}\mathbb{G}_{0}^{v}%
(\theta_{t}^{v}|\psi^{v}),
\end{array}
\end{equation}

Thus $p(\theta_{t}|\mathbf{z}_{1:T},\theta_{-t},\phi)$ can be sampled from with a Metropolis-Hastings (MH) step, where the
candidate pdf is the conditional prior $p(\theta_{t}|\theta_{-t},\phi)$. The
acceptance probability is thus given by
\begin{equation}
\mathsf{\ }\rho(\theta_{t}^{(i)},\theta_{t}^{(i)\ast})=\min\left(
1,\frac{p(\mathbf{z}_{1:T}|\theta_{t}^{(i)\ast},\theta_{-t}^{(i)}%
)}{p(\mathbf{z}_{1:T}|\theta_{t}^{(i)},\theta_{-t}^{(i)})}\right)
\end{equation}
where $\theta_{t}^{(i)\ast}$ is the candidate cluster sampled from
$p(\theta_{t}|\theta_{-t},\phi)$.

The computation of the acceptance probability requires to compute the
likelihood $p(\mathbf{z}_{1:T}|\theta_{t}^{(i)},\theta_{-t}^{(i)})$. This can
be done in $O(T)$ operations using a Kalman filter. However, this has to be
done for $t=1,\ldots,T$ and one finally obtains an algorithm of computational
complexity $O(T^{2})$. Here, we propose to use instead the backward-forward
recursion developed in~\cite{Doucet01}, to obtain an algorithm of overall
complexity $O(T)$. This algorithm uses the following likelihood decomposition
obtained by applying conditional probability rules to $p(\mathbf{z}%
_{1:t-1},\mathbf{z}_{t},\mathbf{z}_{t+1:T}|\theta_{1:T})$
\begin{equation}
p(\mathbf{z}_{1:T}|\theta_{1:T})=p(\mathbf{z}_{1:t-1}|\theta_{1:t-1}%
)p(\mathbf{z}_{t}|\theta_{1:t},\mathbf{z}_{1:t-1})\int_{\mathcal{X}%
}p(\mathbf{z}_{t+1:T}|\mathbf{x}_{t},\theta_{t+1:T})p(\mathbf{x}%
_{t}|\mathbf{z}_{1:t},\theta_{1:t})d\mathbf{x}_{t} \label{eq:backfor}%
\end{equation}
with
\begin{equation}
p(\mathbf{z}_{t:T}|\mathbf{x}_{t-1},\theta_{t:T})\;=\;\int_{\mathcal{X}%
}p(\mathbf{z}_{t+1:T}|\mathbf{x}_{t-1},\theta_{t:T})p(\mathbf{z}%
_{t},\mathbf{x}_{t}|\theta_{t},\mathbf{x}_{t-1})d\mathbf{x}_{t}
\label{eq:backfor2}%
\end{equation}

The first two terms of the r.h.s. in Eq.~\eqref{eq:backfor} are computed by a
forward recursion based on the Kalman filter \cite{Doucet01}. The third term
can be evaluated by a backward recursion according to Eq.~\eqref{eq:backfor2}.
It is shown in \cite{Doucet01} that if $\int_{\mathcal{X}}p(\mathbf{z}%
_{t:T}|\mathbf{x}_{t-1},\theta_{t:T})d\mathbf{x}_{t-1}<\infty$ then
$\frac{p(\mathbf{z}_{t:T}|\mathbf{x}_{t-1},\theta_{t:T})}{\int_{\mathcal{X}%
}p(\mathbf{z}_{t:T}|\mathbf{x}_{t-1},\theta_{t:T})d\mathbf{x}_{t-1}}$ is a
Gaussian distribution w.r.t. $\mathbf{x}_{t-1}$, of mean $m_{t-1|t}^{\prime
}(\theta_{t:T})$ and covariance $P_{t-1|t}^{\prime}(\theta_{t:T})$. Even if $p(\mathbf{z}%
_{t:T}|\mathbf{x}_{t-1},\theta_{t:T})$ is not integrable in $\mathbf{x}_{t-1}$, the
quantities $P_{t-1|t}^{\prime-1}(\theta_{t:T})$ and $P_{t-1|t}^{\prime
-1}(\theta_{t:T})m_{t-1|t}^{\prime}(\theta_{t:T})$ satisfy the backward
information filter recursion (see Appendix). Based on Eq.~\eqref{eq:backfor},
the density $p(\theta_{t}|\mathbf{z}_{1:T},\theta_{-t},\phi)$ is expressed by%
\begin{equation}
p(\theta_{t}|\mathbf{z}_{1:T},\theta_{-t})\;\propto\;p(\theta_{t}|\theta
_{-t},\phi)p(\mathbf{z}_{t}|\theta_{1:t},\mathbf{z}_{1:t-1})\int_{\mathcal{X}%
}p(\mathbf{z}_{t+1:T}|{\mathbf{x}}_{t},\theta_{t+1:T})p(\mathbf{x}%
_{t}|\mathbf{z}_{1:t},\theta_{1:t})d\mathbf{x}_{t}%
\end{equation}

Algorithm 2 summarizes the full posterior sampling procedure. It is the
step-by-step description of Algorithm~1 that accounts for the factorization of
the likelihood given by Eq. (\ref{eq:backfor}).

\bigskip\hrule\smallskip\noindent\textbf{Algorithm 2:} MCMC\ algorithm to
sample from $p(\theta_{1:T}|\mathbf{z}_{1:T},\phi)$\smallskip\hrule\smallskip
\noindent\underline{Initialization $i=1$}

\begin{itemize}
\item For $t=1,...,T$, sample $\theta_{t}^{(1)}.$
\end{itemize}

\noindent\underline{Iteration $i,$ $i=2,\ldots,N^{\prime}+N$}

\begin{itemize}
\item \underline{Backward recursion}: For $t=T,..,1$, compute and store
$P_{t|t+1}^{\prime-1}(\theta_{t+1:T}^{(i-1)})$ and $P_{t|t+1}^{\prime
-1}(\theta_{t+1:T}^{(i-1)})m_{t|t+1}^{\prime}(\theta_{t+1:T}^{(i-1)})$

\item \underline{Forward recursion}: For $t=1,..,T$

\begin{itemize}
\item Perform a Kalman filter step with $\theta_{t}=\theta_{t}^{(i-1)}$, store
$\widehat{\mathbf{x}}_{t|t}(\theta_{1:t-1}^{(i)},\theta_{t}^{(i-1)})$ and
$\Sigma_{t|t}(\theta_{1:t-1}^{(i)},\theta_{t}^{(i-1)}).$

\item Metropolis-Hastings step :

\begin{itemize}
\item Sample a candidate cluster%
\begin{equation}
\theta_{t}^{(i)\ast}\sim p(\theta_{t}|\theta_{-t}^{(i)},\phi)
\end{equation}

\item Perform a Kalman filter step with $\theta_{t}=\theta_{t}^{(i)\ast}$,
store $\widehat{\mathbf{x}}_{t|t}(\theta_{1:t-1}^{(i)},\theta_{t}^{(i)\ast})$
and $\Sigma_{t|t}(\theta_{1:t-1}^{(i)},\theta_{t}^{(i)\ast})$

\item Compute
\begin{equation}
\rho(\theta_{t}^{(i)},\theta_{t}^{(i)\ast})=\min\left(  1,\frac{p(\mathbf{z}%
_{1:T}|\theta_{t}^{(i)\ast},\theta_{-t}^{(i)})}{p(\mathbf{z}_{1:T}|\theta
_{t}^{(i)},\theta_{-t}^{(i)})}\right)
\end{equation}

\item With probability $\rho(\theta_{t}^{(i)},\theta_{t}^{(i)\ast})$, set
$\theta_{t}^{(i)}=\theta_{t}^{(i)\ast}$, otherwise $\theta_{t}^{(i)}%
=\theta_{t}^{(i-1)}.$
\end{itemize}
\end{itemize}
\end{itemize}

\noindent\underline{State post-Sampling (for non-burn-in iterations only)}

\begin{itemize}
\item For $i=N^{\prime}+1,...,N^{\prime}+N$, compute $\widehat{\mathbf{x}%
}_{t|T}(\theta_{1:T}^{(i)})=\mathbb{E}\left(  \mathbf{x}_{t}|\theta
_{1:T}^{(i)},\mathbf{z}_{1:T}\right)  $ for all $t$ with a Kalman smoother.
\end{itemize}


\smallskip\hrule\bigskip

It can be easily established that the simulated Markov chain $\left\{
\theta_{1:T}^{\left(  i\right)  }\right\}  $ is ergodic with limiting
distribution $p(\theta_{1:T}|\mathbf{z}_{1:T})$. After $N^{\prime}$ burn-in,
the $N$ last iterations of the algorithm are kept, and the MMSE estimates of
$\theta_{t}$ and $\mathbf{x}_{t}$ for all $t=0,\ldots,T$ are computed as
explained in Subsection~\ref{sec:mixt-param-estim}, using
\begin{equation}
\widehat{\theta}_{t|T}^{\text{{\tiny MMSE}}}=\frac{1}{N}\sum_{i=N'+1}^{N'+N}%
\theta_{t}^{(i)}\quad\widehat{\mathbf{x}}_{t|T}^{\text{{\tiny MMSE}}}=\frac
{1}{N}\sum_{i=N'+1}^{N'+N}\widehat{\mathbf{x}}_{t|T}(\theta_{1:T}^{(i)})
\end{equation}

\subsection{Unknown hyperparameters}

\label{sec:unkn-hyperp} The hyperparameters in vector $\phi$ have some
influence on the correct estimation of the DPMs $F^{v}$ and $F^{w}$. In this
subsection, we include them in the inference by considering them as unknowns
with prior distributions:
\begin{align}
\alpha^{v}  &  \sim\mathcal{G}(\frac{\eta}{2},\frac{\nu}{2}),\ \alpha^{w}%
\sim\mathcal{G}(\frac{\eta}{2},\frac{\nu}{2}),\ \\
\psi^{v}  &  \sim p_{0}(\psi^{v}),\ \psi^{w}\sim p_{0}(\psi^{w})
\end{align}
where $\eta$ and $\nu$ are known constants and $p_{0}$ is a pdf with fixed and
known parameters. The posterior probability $p(\alpha^{v}|\mathbf{x}%
_{1:T},\theta_{1:T},\mathbf{z}_{1:T},\psi^{v},\phi^{w})$ reduces to
$p(\alpha^{v}|M^{v},T)$ where $M^{v}$ is the number of distinct values taken
by the clusters $\theta_{1:T}^{v}$. As shown in~\cite{Antoniak74}, this pdf
can be expressed by
\begin{equation}
p(\alpha^{v}|M^{v},T)\propto\frac{s(T,M^{v})(\alpha^{v})^{M^{v}}}{\sum
_{k=1}^{T}s(T,k)(\alpha^{v})^{k}}p(\alpha^{v}) \label{eq:postalpha}%
\end{equation}
where the $s(T,k)$ are the absolute values of Stirling numbers of the first
kind. We can sample from the above pdf with a Metropolis-Hasting step using
the prior Gamma pdf $p(\alpha^{v})=\mathcal{G}(\frac{\eta}{2},\frac{\nu}{2})$
as proposal (and similarly for $\alpha^{w}$). Other methods have been proposed
that allow direct sampling, see for example West \cite{West92},\ and Escobar
and West \cite{Escobar95}.

The posterior probability $p(\psi^{v}|\mathbf{x}_{1:T},\theta_{1:T}%
,\mathbf{z}_{1:T},\alpha^{v},\phi^{w})$ reduces to $p(\psi^{v}|\theta
_{1:M^{v}}^{v\ \prime})$ where $\theta_{1:M^{v}}^{v\ \prime}$ is the set of
distinct values taken by the clusters $\theta_{1:T}^{v}$. It is expressed by
\begin{equation}
p(\psi^{v}|\mathbf{x}_{1:T},\theta_{1:T},\mathbf{z}_{1:T},\alpha^{v},\phi
^{w})\propto p_{0}(\psi^{v})\prod_{k=1}^{M^{v}}\mathbb{G}_{0}^{v}(\theta
_{k}^{v\ \prime}|\psi^{v}) \label{eq:postpsi}%
\end{equation}
We can sample from this pdf with a Metropolis-Hasting step using the prior
Gamma pdf $p_{0}(\psi^{v})$ as proposal whenever direct sampling is not possible.

\section{Rao-Blackwellized Particle Filter algorithm for online state
estimation}

\label{sec:SMC} Many applications, such as target tracking, require
\emph{online} state estimation. In this case, the MCMC\ approach is inadequate
as it requires availability of the entire dataset to perform state estimation.
In this section, we develop the online counterpart to the MCMC procedure
presented in Section~\ref{sec:MCMC}:  a sequential Monte Carlo
method (also known as particle filter) is implemented, to sample on-line from
the sequence of probability distributions $\{ p(\mathbf{x}_{0:t},\theta
_{1:t}|\mathbf{z} _{1:t})$, $t=1,2,\ldots\} $. Here, the
hyperparameter vector $\phi$ is assumed to be known, therefore it is omitted
in the following. Online hyperparameter estimation is discussed in
Section~\ref{sec:discussion}.

As explained in Subsection~\ref{sec:mixt-param-estim}, we need to sample from
$p(\theta_{1:t}|\mathbf{z}_{1:t})$, because $p(\mathbf{x}_{0:t}|\theta
_{1:t},\mathbf{z}_{1:t})$ can be computed using Kalman techniques. (The
sampling procedure is indeed a generalization of the Rao-Blackwellized
particle filter~\cite{Doucet01b} to DPMs.) At time $t$, $p(\mathbf{x}%
_{t},\theta_{1:t}|\mathbf{z}_{1:t})$ is approximated through a set of $N$
particles $\theta_{1:t}^{(1)},\ldots,\theta_{1:t}^{(N)}$ by the\ following
empirical distribution%
\begin{equation}
P_{N}(\mathbf{x}_{t},\theta_{1:t}|\mathbf{z}_{1:t})=\sum_{i=1}^{N}{w}%
_{t}^{(i)}\mathcal{N}(\mathbf{x}_{t};\widehat{\mathbf{x}}_{t|t}(\theta
_{1:t}^{(i)}),\Sigma_{t|t}(\theta_{1:t}^{(i)})) \label{eq:empiricaldistrib}%
\end{equation}
The parameters $\widehat{\mathbf{x}}_{t|t}(\theta_{1:t}^{(i)})$ and
$\Sigma_{t|t}(\theta_{1:t}^{(i)})$ are computed recursively for each particle
$i$ using the Kalman filter \cite{Anderson79}. In order to build the
algorithm, we note that%
\begin{equation}
p(\theta_{1:t}^{(i)}|\mathbf{z}_{1:t})\propto p(\theta_{1:t-1}^{(i)}%
|\mathbf{z}_{1:t-1})p(\mathbf{z}_{t}|\theta_{1:t}^{(i)},\mathbf{z}%
_{1:t-1})p(\theta_{t}^{(i)}|\theta_{1:t-1}^{(i)})
\end{equation}
where
\begin{equation}%
\begin{array}
[c]{lcl}%
p(\mathbf{z}_{t}|\theta_{1:t}^{(i)},\mathbf{z}_{1:t-1}) & = & p(\mathbf{z}%
_{t}|\theta_{t}^{(i)},\theta_{1:t-1}^{(i)},\mathbf{z}_{1:t-1})\nonumber\\
& = & \mathcal{N}(\mathbf{z}_{t};\widehat{\mathbf{z}}_{t|t-1}(\theta
_{1:t}^{(i)}),S_{t|t-1}(\theta_{1:t}^{(i)}))
\end{array}
\end{equation}

and%
\begin{align}
\widehat{\mathbf{z}}_{t|t-1}(\theta_{1:t}^{(i)})  &  =H_{t}\left[  F_{t} \,
\widehat{\mathbf{x}}_{t-1|t-1}(\theta_{1:t-1}^{(i)}) +C_{t}\mathbf{u}%
_{t}+G_{t}\mu_{t}^{v\ (i)}\right]  +\mu_{t}^{w\ (i)}\\
S_{t|t-1}(\theta_{1:t}^{(i)})  &  =H_{t}\left[  F_{t}\, \Sigma_{t-1|t-1}%
(\theta_{1:t-1}^{(i)})\, F_{t}^{\mathsf{T}}+G_{t}\Sigma_{t}^{v\ (i)}%
G_{t}^{\mathsf{T}}\right]  H_{t}^{\mathsf{T}}+\Sigma_{t}^{w\ (i)}\nonumber
\end{align}

The Rao-Blackwellized Particle Filter (RBPF) algorithm proceeds as follows.

\bigskip\hrule\smallskip

\noindent\textbf{Algorithm 3:} Rao-Blackwellized Particle Filter to sample
from $p(\theta_{1:t}|\mathbf{z}_{1:t})$ \smallskip\hrule\smallskip


\underline{At time $0$}$.$

\begin{itemize}
\item For $i=1,..,N$, sample $\left(  \widehat{\mathbf{x}}_{0|0}^{(i)}%
,\Sigma_{0|0}^{(i)}\right)  \sim p_{0}(\mathbf{x}_{0|0},\Sigma_{0|0}).$

\item Set $w_{0}^{(i)}\leftarrow\frac{1}{N}$
\end{itemize}

\underline{At each time $t$} ($t\geq1$), do for $i=1,\ldots,N$

\begin{itemize}
\item Sample $\widetilde{\theta}_{t}^{(i)}\sim q(\theta_{t}|\theta
_{1:t-1}^{(i)},\mathbf{z}_{1:t})$

\item Compute $\{\widehat{\mathbf{x}}_{t|t-1}(\theta_{1:t-1}^{(i)},\widetilde{\theta}_{t}^{(i)}
),\Sigma_{t|t-1}(\theta_{1:t-1}^{(i)},\widetilde{\theta}_{t}^{(i)}),\widehat{\mathbf{x}}_{t|t}(\theta_{1:t-1}^{(i)},\widetilde{\theta}_{t}^{(i)}),\Sigma_{t|t}(\theta_{1:t-1}^{(i)},\widetilde{\theta}_{t}^{(i)})\}$ by using a Kalman filter
step from $\{\widehat{\mathbf{x}}_{t-1|t-1}(\theta_{1:t-1}^{(i)})$,
$\Sigma_{t-1|t-1}(\theta_{1:t-1}^{(i)})$, $\widetilde{\theta}_{t}^{(i)}$,
$\mathbf{z}_{t})\}$

\item For $i=1,\ldots,N$, update the weights according to
\begin{equation}
\widetilde{w}_{t}^{(i)}\varpropto w_{t-1}^{(i)}\frac{p(\mathbf{z}_{t}%
|\theta_{1:t-1}^{(i)},\widetilde{\theta}_{t}^{(i)},\mathbf{z}_{1:t-1}%
)p(\widetilde{\theta}_{t}^{(i)}|\theta_{1:t-1}^{(i)})}{q(\widetilde{\theta
}_{t}^{(i)}|\theta_{1:t-1}^{(i)},\mathbf{z}_{1:t})}%
\end{equation}

\item Compute $S=\sum_{i=1}^{N}\widetilde{w}_{t}^{(i)}$ and for $i=1,\ldots
,N$, set $\widetilde{w}_{t}^{(i)}\leftarrow\frac{\widetilde{w}_{t}^{(i)}}{S}$

\item Compute $N_{\text{eff}}=\left[  \sum_{i=1}^{N}\left(  \widetilde{w}%
_{t}^{(i)}\right)  ^{2}\right]  ^{-1}$

\item If $N_{\text{eff}}\leq\eta$, then resample the particles -- that is,
duplicate the particles with large weights are remove the particles with small
weights. This results in a new set of particles denoted ${\theta}_{t}^{(i)}$
with weights $w_{t}^{(i)}=\frac{1}{N}$

\item Otherwise, rename the particles and weights by removing the
$\widetilde{\cdot}$'s.
\end{itemize}

\smallskip\hrule\bigskip

Particle filtering convergence results indicate that the variance of the Monte Carlo estimates depends highly on the importance distribution selected. Here, the conditionally optimal importance distribution is $q(\theta_{t}|\theta_{1:t-1}^{(i)}%
,\mathbf{z}_{1:t})=p(\theta_{t}|\theta_{1:t-1}^{(i)},\mathbf{z}_{1:t})$,
see~\cite{Doucet01b}. However, it cannot be used, as the associated importance
weights do not admit a closed-form expression\footnote{When using the optimal
importance distribution, the weights computation requires the evaluation of an
integral with respect to $\theta_{t}$. It is possible to integrate
analytically w.r.t. the cluster means $\mu^{v}$ and $\mu^{w}$, but not w.r.t.
the covariances.}. In practice, the evolution pdf $p(\theta_{t}|\theta
_{1:t-1})$ was used as the importance distribution.

From the particles, the MMSE estimate and posterior covariance matrix of
$\mathbf{x}_{t}$ are given by%
\begin{equation}
\widehat{\mathbf{x}}_{t|t}^{\text{{\tiny MMSE}}}=\sum_{i=1}^{N}{w}_{t}%
^{(i)}\widehat{\mathbf{x}}_{t|t}(\theta_{1:t}^{(i)})
\end{equation}%
\begin{equation}
\widehat{\Sigma}_{t|t}=\sum_{i=1}^{N}{w}_{t}^{(i)}\left[  \Sigma_{t|t}%
(\theta_{1:t}^{(i)})+(\widehat{\mathbf{x}}_{t|t}(\theta_{1:t}^{(i)}%
)-\widehat{\mathbf{x}}_{t|t}^{\text{{\tiny MMSE}}})(\widehat{\mathbf{x}}%
_{t|t}(\theta_{1:t}^{(i)})-\widehat{\mathbf{x}}_{t|t}^{\text{{\tiny MMSE}}%
})^{T}\right]
\end{equation}

\section{Applications}

\label{sec:applications} In this section, we present two applications of the
above model and algorithms\footnote{See Caron et al.~\cite{Caron2006} for an application on a regression problem.}. We address, first, blind deconvolution, second,
change point detection in biomedical time series. In each case, we assume that
the statistics of the state noise are unknown, and modelled as a DPM.

\subsection{Blind deconvolution of impulse processes}

Various fields of Engineering and Physics, such as image de-blurring,
spectroscopic data analysis, audio source restoration, etc. require blind
deconvolution. We follow here the model presented in~\cite{Doucet97} for blind deconvolution of Bernoulli-Gaussian processes, which is
recalled below.

\subsubsection{Statistical Model}

Let $H=\left(
\begin{array}
[c]{cccc}%
1 & h_{1} & .. & h_{L}%
\end{array}
\right)  =\left(
\begin{array}
[c]{cc}%
1 & \mathbf{h}%
\end{array}
\right)  $ and $\mathbf{x}_{t}=\left(
\begin{array}
[c]{cccc}%
v_{t} & v_{t-1} & ... & v_{t-L}%
\end{array}
\right)  ^{T}$. The observed signal $z_{t}$ is the convolution of the sequence
$\mathbf{x}_{t}$ with a finite impulse response filter $H$, observed in
additive white Gaussian noise $w_{t}$. The observation model is then%
\begin{equation}
z_{t}=H\mathbf{x}_{t}+w_{t}%
\end{equation}
where $w_{t}\sim\mathcal{N}(0,\sigma_{w}^{2})$ with $\sigma_{w}^{2}$ is the
assumed known variance of $w_{t}$. The state space model can be written as
follows:
\begin{equation}
\mathbf{x}_{t}=F\mathbf{x}_{t-1}+Gv_{t}%
\end{equation}
where $F=\left(
\begin{array}
[c]{cc}%
0 & 0_{1\times L}\\
0_{L\times1} & I_{L}%
\end{array}
\right)  $, $G=\left(
\begin{array}
[c]{c}%
1\\
0_{L\times1}%
\end{array}
\right)  $, $0_{m\times n}$ is the zero matrix of size $m\times n$\ and
$I_{m}$ is the identity matrix of size $m\times m$. The state transition noise
$v_{t}$ is supposed to be independent from $w_{t}$, and distributed according
to the mixture%
\begin{equation}
v_{t}\;\sim\;\lambda F^{v}+(1-\lambda)\delta_{0} \label{eq:blindstuffvt}%
\end{equation}
where $\delta_{0}$ is the Dirac delta function at $0$ and $F^{v}$ is a DPM of
Gaussians defined in Eq.~\eqref{eq:defGaussianDPM}. In other words, the noise
is alternatively zero, or distributed according to a DPM of Gaussians.

For simplicity reasons, we introduce latent Bernoulli variables $r_{t}%
\in\{0,1\}$ such that $\Pr(r_{t}=1)=\lambda$ and $v_{t}|(r_{t}=1)\sim
f(\cdot|\theta_{t}^{v}),\ v_{t}|(r_{t}=0)\sim\delta_{0}$. Consider the cluster
variable $\varphi_{t}^{v}$ defined by $\varphi_{t}^{v}=\theta_{t}^{v}$ if $r_{t}=1$
and $\varphi_{t}^{v}=(0,0)$ (i.e. parameters corresponding to the delta-mass) if
$r_{t}=0$, that is, $\varphi_{t}^{v}\sim\lambda F^{v}+(1-\lambda)\delta_{(0,0)}$.
By integrating out $F^{v}$, one has%
\begin{equation}
\varphi_{t}^{v}|\varphi_{-t}^{v}\sim\lambda p(\varphi_{t}^{v}|\varphi_{-t}^{v}%
,r_{t}=1)+(1-\lambda)\delta_{(0,0)} \label{eq:blindstuffphi1}%
\end{equation}
where $p(\varphi_{t}^{v}|\varphi_{-t}^{v},r_{t}=1)$ is the Polya urn representation
on the set $\widetilde{\varphi}_{-t}^{v}=\{\varphi\in\varphi_{-t}^{v}|\varphi\neq
\delta_{(0,0)}\}$ of size $T^{\prime}$ given by
\begin{equation}
\varphi_{t}^{v}|(\varphi_{-t}^{v},r_{t}=1)\sim\frac{\sum_{k=1,k\neq t}^{T^{\prime}%
}\delta_{\varphi_{k}^{v}}+\alpha^{v}\mathbb{G}_{0}^{v}}{\alpha
^{v}+T^{\prime}}%
\end{equation}
\bigskip

The probability $\lambda$ is considered as a random variable with a beta prior
density $p(\lambda)=\mathcal{B}(\zeta,\tau)$ where $\zeta$ and $\tau$ are
known parameters. The random variable $\lambda$ can be marginalized out in
Eq.~\eqref{eq:blindstuffphi1}%
\begin{equation}
\varphi_{t}^{v}|\varphi_{-t}^{v}\sim\frac{a(\varphi_{-t}^{v})}{a(\varphi_{-t}^{v}%
)+b(\varphi_{-t}^{v})}p(\varphi_{t}^{v}|\varphi_{-t}^{v},r_{t}=1)+\frac{b(\varphi_{-t}%
^{v})}{a(\varphi_{-t}^{v})+b(\varphi_{-t}^{v})}\delta_{(0,0)}%
\end{equation}
where%
\begin{align}
a(\varphi_{-t}^{v})  &  =\zeta+\sum_{k=1,k\neq t}^{T}r_{k}\\
b(\varphi_{-t}^{v})  &  =\tau+\sum_{k=1,k\neq t}^{T}(1-r_{k})
\end{align}
where $r_{t}=0$ if $\varphi_{t}^{v}=(0,0)$ and $r_{t}=1$ otherwise.

The hyperparameters are $\phi=(\alpha^{v},\mathbf{h})$ (the hyperparameters of
the base distribution $\mathbb{G}_{0}^{v}$ are assumed fixed and known). These
hyperparameters are assumed random with prior distribution $p(\phi
)=p(\alpha^{v})p(\mathbf{h})$, where
\begin{equation}
p(\alpha^{v})=\mathcal{G}(\frac{\eta}{2},\frac{\nu}{2}),\ p(\mathbf{h}%
)=\mathcal{N}(0,\sigma_{w}^{2}\Sigma_{\mathbf{h}})
\end{equation}
where $\eta$, $\nu$ and $\Sigma_{\mathbf{h}}$ are known. Conditional on
$\mathbf{x}_{0:t}$, the following conditional posterior is obtained
straighforwardly%
\begin{equation}
p(\mathbf{h}|\mathbf{x}_{0:t},\mathbf{z}_{1:T})=\mathcal{N}(\mathbf{m}%
,\sigma_{w}^{2}\Sigma_{\mathbf{h}}^{\prime})
\end{equation}
where
\begin{align*}
\Sigma_{\mathbf{h}}^{\prime-1}  &  =\Sigma_{\mathbf{h}}^{-1}+\sum_{t=1}%
^{T}\mathbf{v}_{t-1:t-L}\mathbf{v}_{t-1:t-L}^{\prime}\\
\mathbf{m}  &  =\Sigma_{\mathbf{h}}^{\prime}\sum_{t=1}^{T}\mathbf{v}%
_{t-1:t-L}\left(  z_{t}-v_{t}\right)
\end{align*}

Samples $\mathbf{x}_{0:t}^{(i)}$ can be generated from the Gaussian posterior
$p(\mathbf{x}_{0:t}|\varphi_{1:T}^{v\ (i)},\mathbf{z}_{1:T},\phi^{(i-1)})$ with
the simulation smoother~\cite{Durbin2002}. This algorithm complexity is
$O(T)$.

The aim is to approximate by MCMC the joint posterior pdf $p(\mathbf{v}%
_{1:T},\varphi_{1:T},\phi|\mathbf{z}_{1:T})$. This is done by implementing
Algorithm~3 for the cluster variable, whereas the other variables are sampled
by Metropolis-Hastings or direct sampling w.r.t their conditional posterior.

\subsubsection{Simulation results}

\begin{figure}
[h]
\begin{center}
\includegraphics[width=3.5in]%
{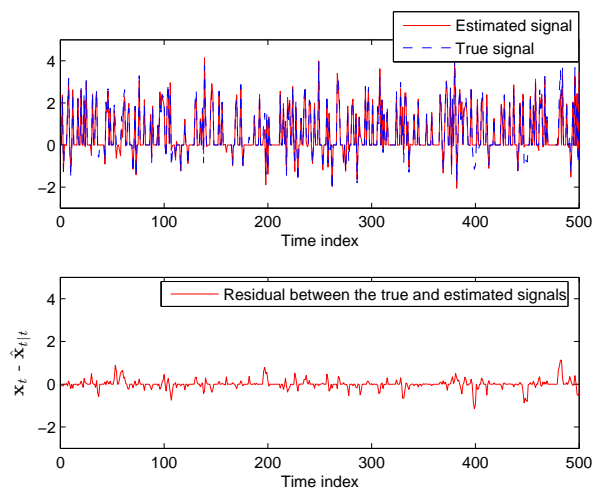}%
\caption{Top picture: True (dashed line) and MMSE estimated (solid line)
signal $\mathbf{v}_{1:T}$ after 10,000 MCMC iterations (7,500 burn-in).
$v_{t}$ is supposed to be either $0$ with probability $\lambda$, or to be
distributed from an unknown pdf $F^{v}$ with probability $(1-\lambda)$. Bottom
picture: residual $e_{t}=v_{t}-E[v_{t}|\mathbf{z}_{1:T}]$ between the true and
estimated signals. Although the distribution $F^{v}$ is unknown, the state
$v_{t}$ is almost correctly estimated.}%
\label{fig:deconv1}%
\end{center}
\end{figure}

\begin{figure}
[hh]
\begin{center}
\includegraphics[width=3.5in]%
{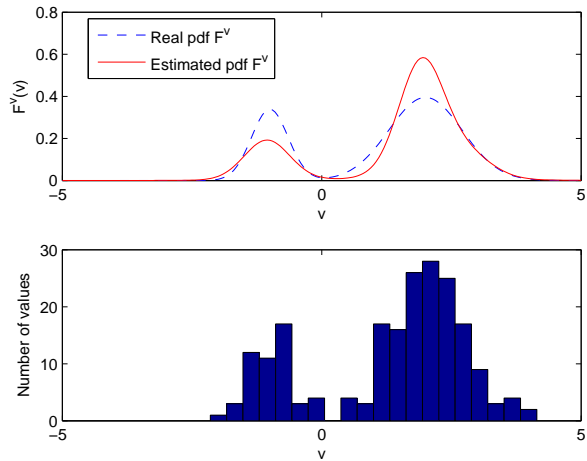}%
\caption{(Top) True (dashed line) and estimated (solid line) pdf $F^{v}$. The
true pdf $F^{v}$ is a mixture of two Gaussians $0.7\mathcal{N}%
(2,.5)+0.3\mathcal{N}(-1,.1)$. It is supposed to be unknown and jointly
estimated with the state vector with $10,000$ MCMC iterations (7,500 burn-in)
given a vector of $120$ observations $\mathbf{z}_{1:T}$. The estimated pdf
matches correctly the two modes of the true distribution. (Bottom) Histogram
of the simulated values $v$, sampled from $F^{v}$ which a mixture of two Gaussians
$0.7\mathcal{N}(2,.5)+0.3\mathcal{N}(-1,.1)$}%
\label{fig:deconv2}%
\end{center}
\end{figure}

This model has been simulated with the following parameters: $T=120$, $L=3$,
$\mathbf{h}=\left(
\begin{array}
[c]{ccc}%
-1.5 & 0.5 & -0.2
\end{array}
\right)  $, $\lambda=0.4$, $\sigma_{w}^{2}=0.1,F^{v}=0.7\mathcal{N}%
(2,.5)+0.3\mathcal{N}(-1,.1)$, $\Sigma_{\mathbf{h}}=100$, $\eta=3$, $\nu=3$,
$\zeta=1$, $\tau=1$. The hyperparameters of the base distribution are $\mu
_{0}=0$, $\kappa_{0}=0.1$, $\nu_{0}=4,\Lambda_{0}=1$. For the estimation,
10,000 MCMC iterations are performed, with 7,500 burn-in iterations.
Fig.~\ref{fig:deconv1} (top) displays the MMSE estimate of $\mathbf{v}_{1:T}$
together with its true value. As can be seen in Fig.~\ref{fig:deconv1}
(bottom), the signal is correctly estimated and the residual is quite small.
Also, as can be seen in Fig.~\ref{fig:deconv2}, the estimated pdf $F^{v}$ is
quite close to the true one. In particular, the estimated pdf matches the two
modes of the true pdf. Multiple simulations with different starting values were runned, and the results appeared insensitive to initialization. This suggest that the MCMC sampler explores properly the posterior.%
\begin{figure}
[hhh]
\begin{center}
\includegraphics[
height=2.5036in,
width=3.3243in
]%
{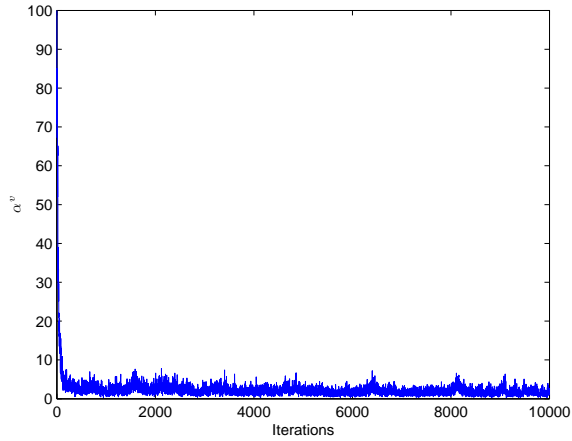}%
\caption{Evolution of $\alpha^{v\ (i)}$ in function of Gibbs sampler iteration
i. The value of $\alpha^{v}$ is initialized at 100.}%
\label{fig:deconv3}%
\end{center}
\end{figure}
\begin{figure}
[hhhh]
\begin{center}
\includegraphics[
height=2.5685in,
width=3.4117in
]%
{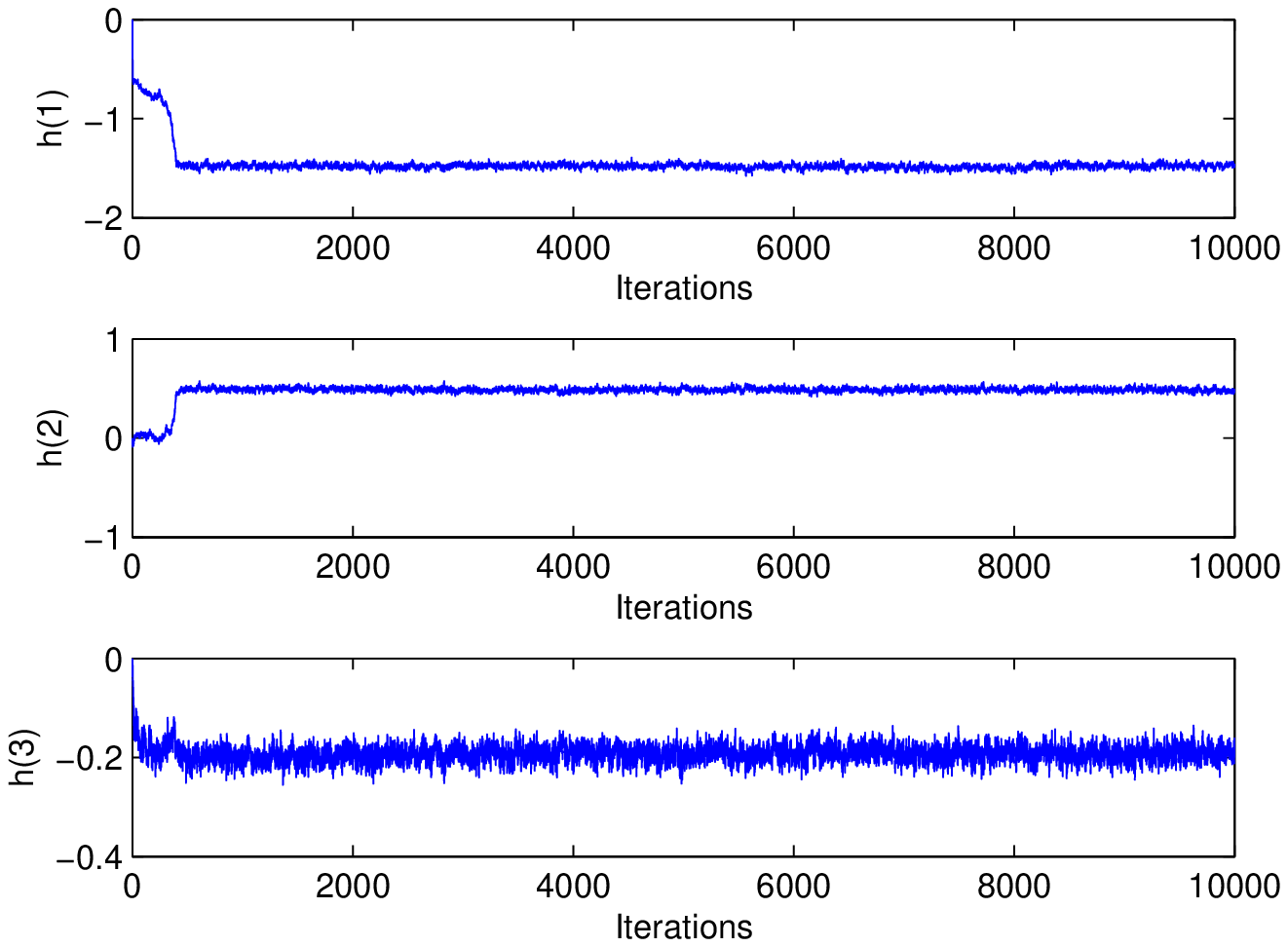}%
\caption{Evolution of the three components of the vector $\mathbf{h}^{(i)}$ in
function of Gibbs sampler iteration i. It is initialized at $[0~0~0]$. The
value converges toward the true value $\mathbf{h}=[-1.5\ 0.5\ -0.2]$.}%
\label{fig:deconv4}%
\end{center}
\end{figure}

Let $e_{MSE}$ be the mean squared error (MSE), computed by%
\begin{equation}
e_{MSE}=\sqrt{\frac{1}{T}\sum_{t=1}^{T}(\mathbf{v}_{t}-\mathbf{v}%
_{t|T}^{\text{{\footnotesize MMSE}}})^{2}}%
\end{equation}

To better highlight the performance of the proposed algorithm, we compared our
model/algorithm (denoted M1) with the following models, denoted M2 to M8:

\begin{description}
\item[M2.] In this model, the pdf is assumed known and set to the true value
$F^{v}=0.7\mathcal{N}(2,.5)+0.3\mathcal{N}(-1,.1)$. The model is simply a Jump
Linear Model that jumps between three modes of resp. mean/covariance $(0,0)$,
$(2,.5)$ and $(-1,.1)$ with resp. prior probabilities $(1-\lambda)$,
$0.7\lambda$ and $0.3\lambda$.

\item[M3.] In this model, the pdf is assumed to be a Gaussian
$\mathcal{N}(1.1,2.3)$. The first two moments of this Gaussian are the same as
those of the true pdf $F^{v}$. The model is also a Jump Linear Model that
jumps between two modes\ of resp. mean/covariance $(0,0)$ and $(1.1,2.3)$ with
resp. prior probabilities $(1-\lambda)$ and $\lambda$.

\item[M4-7.] The model described in this article but with $\alpha^{v}$ fixed to
$0.1$ (M3), $1$ (M4), $10$ (M5) and $100$ (M6).

\item[M8.] The model described in this article (M1) but with the observation
noise variance $\sigma_{w}^{2}$ estimated with an inverse gamma prior
$\sigma_{w}^{2}\sim i\mathcal{G}(u,v)$ with $u=2$ and $v=0.1$. $\sigma
_{w}^{2\ (i)}$ is sampled with Gibbs sampling with $\sigma_{w}^{2}%
|\mathbf{x}_{0:T},\mathbf{z}_{1:T},\mathbf{h\sim}i\mathcal{G}(u^{\prime
},v^{\prime})$ and $u^{\prime}=u+\frac{T}{2}$ and $v^{\prime}=v+\frac{1}%
{2}\sum_{t=1}^{T}(\mathbf{z}_{t}-H\mathbf{x}_{t})^{2}$.
\end{description}

The algorithm used for M2 and M3 is the Gibbs sampler with backward forward
recursion given in~\cite{Doucet01}. For the same set of observations, each MCMC algorithm has been run with
10,000 iterations and 7,500 burn-in iterations. MMSE estimate $\mathbf{v}%
_{t|T}^{\text{{\footnotesize MMSE}}}$ and MSE $e_{MSE}$ are computed for each
model. 20 simulations have been performed; for each model, the mean and standard deviation of the MSE's over the 20 simulations are reported in Tab. I.

\begin{center}
Tab. I. Comparison of our model/algorithm with other models%

\begin{tabular}
[c]{|l|l|l|l|l|l|l|l|l|}\hline
Simulation / Model & M1 & M2 & M3 & M4 & M5 & M6 & M7 & M8\\\hline\hline
Mean & 0.240&   0.217&  0.290&  0.915&  0.254   &0.253  &0.314  &0.438
\\\hline
Standard deviation & 0.067& 0.058&  0.085&  0.818&  0.062&  0.086&  0.222&  0.421
\\\hline
\end{tabular}

\end{center}

Our model/algorithm (M1) gives MSE that is only 10\% more than that of
the model with fixed pdf (M2) even though the pdf is not exactly
estimated. If the observation noise variance $\sigma_{w}^{2}$ is
unknown and has to be estimated (M8), this has an impact on the
estimation of the state vector still the sampler converge more slowly
to the true posterior. If the unknown pdf is set to be a Gaussian with
large variance (M3), the MSE is 17\% larger than with our approach.
The estimation of $\alpha^{v}$ improves the estimation of the state
vector: MSEs are higher for models M4-7 where $\alpha^{v}$ is set to a
fixed value. This is especially true for $\alpha^{v}=0.1$. With this
small value, the sampler proposes new clusters very rarely and
converges very slowly to the true posterior.

\subsection{Change-point problems in biomedical time series}

Let now consider a change-point problem in biomedical time series. The
following problem has been discussed in~\cite{Gordon90} and~\cite{Carter96}.
Let consider patients who had recently undergone kidney transplant. The level
of kidney function is given by the rate at which chemical substances are
cleared from the blood, and the rate can be inferred indirectly from
measurements on serum creatinine. If the kidney function is stable, the
response series varies about a constant level. If the kidney function is
improving (resp. decaying) at a constant level then the response series decays
(resp. increases) linearly.

\subsubsection{Statistical model}

The linear model, formulated by Gordon and Smith~\cite{Gordon90} is given by
\begin{align}
\mathbf{x}_{t}  &  =F\mathbf{x}_{t-1}+G\mathbf{v}_{t}\\
z_{t}  &  =H\mathbf{x}_{t}+w_{t}%
\end{align}
where $\mathbf{x}_{t}=(m_{t},\dot{m}_{t}),$ where $m_{t}$ is the level and
$\dot{m}_{t}$ the slope, $F=\left(
\begin{array}
[c]{cc}%
1 & 1\\
0 & 1
\end{array}
\right)  $, $G=\left(
\begin{array}
[c]{cc}%
1 & 1\\
0 & 1
\end{array}
\right)  $, $z_{t}$ is the measured creatinine and $H=\left(
\begin{array}
[c]{cc}%
1 & 0
\end{array}
\right)  $. Measurements are subject to errors due to mistakes in data
transcription, equipment malfunction or blood contamination. $w_{t}$ follows
the following mixture model
\begin{equation}
w_{t}\sim\lambda^{w}\mathcal{N}(0,\sigma_{1}^{w})+(1-\lambda^{w}%
)\mathcal{N}(0,\sigma_{2}^{w})
\end{equation}
where $\lambda^{w}=0.98$ is the probability that the measurements are correct,
in that case the variance is $\sigma_{1}^{w}=10^{-7}$ and $\sigma_{2}^{w}=1$
otherwise. To capture the effects of jumps in the creatinine level, the state
noise $\mathbf{v}_{t}$ is supposed to be distributed according to the
following mixture model%
\begin{equation}
\mathbf{v}_{t}\sim\lambda^{v}F^{v}+(1-\lambda^{v})\delta_{\theta_{0}^{v}}%
\end{equation}
where $\theta_{0}^{v}=\left\{  \left(
\begin{array}
[c]{cc}%
0 & 0
\end{array}
\right)  ^{T},\left(
\begin{array}
[c]{cc}%
0 & 0\\
0 & 0
\end{array}
\right)  \right\}  $, $\lambda^{v}=0.15$ is the probability of jump in the
level and $F^{v}$ is a DPM of Gaussians. Contrary to the model
in~\cite{Carter96}, we do not define fixed jump levels. These levels, as well
as their number, are estimated through the DPM.

\subsubsection{Simulation results}

\begin{figure}
[ptb]
\begin{center}
\includegraphics[
height=2.0064in,
width=2.6662in
]%
{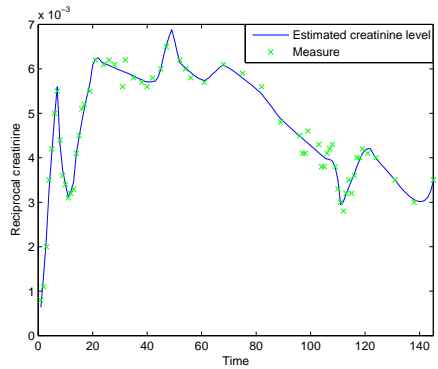}%
\caption{Measured (cross) and estimated (solid line) creatinine level with
2000 MCMC iterations and 1000 burn-in iterations.}%
\label{fig:biomed1}%
\end{center}
\end{figure}

\begin{figure}
[ptbptb]
\begin{center}
\includegraphics[
height=1.9752in,
width=2.6212in
]%
{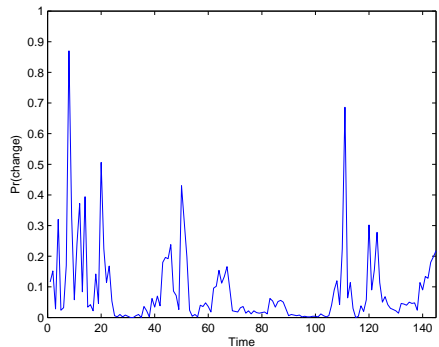}%
\caption{Posterior probability of a jump in the creatinine level with 2000
MCMC iterations and 1000 burn-in iterations. For a threshold set to $0.5$, the creatinine level experiences jumps at about times $8$, $20$ and $110$.}%
\label{fig:biomed2}%
\end{center}
\end{figure}

The last model is applied to the data provided in Gordon and Smith
\cite{Gordon90} (and also exploited in~\cite{Carter96}). The hyperparameters
of the base distribution are $\mu_{0}=\left(
\begin{array}
[c]{c}%
0\\
0
\end{array}
\right)  $, $\kappa_{0}=10^{6}$, $\nu_{0}=4$, $\Lambda_{0}=\frac{10^{-6}}%
{2}\left(
\begin{array}
[c]{cc}%
1 & 0\\
0 & 1
\end{array}
\right)  $. For the estimation, 2,000 MCMC iterations (with 1,000 burn-in
iterations) are performed. Fig.~\ref{fig:biomed1}\ presents the estimated
creatinine level together with the measurements. Fig.~\ref{fig:biomed2}\ plots
the posterior probability of a jump in the creatinine level. In particular, the estimated pdf matches the two
modes of the true pdf. Multiple simulations with different starting values were runned, and the results appeared insensitive to initialization. This suggest that the MCMC sampler explores properly the posterior.%

The estimation have also been made online with the
Rao-Blackwellized algorithm with 1000 particles. We perform
fixed-lag smoothing~\cite{Doucet01c} to estimate
$\mathbb{E}(\mathbf{x}_t|\mathbf{z}_{1:t+T})$, where $T$ is set to
$10$. The mean time per iteration is about 1s. The importance
function used to sample the latent variables $\theta_{t}^{v}$ is
prior pdf $p(\theta_{t}^{v}|\theta_{1:t-1}^{v})$. For a detection
threshold set at $0.5$, the MCMC algorithm detects 3 peaks, while
the RBPF only detects two peaks. The trade-off between false alarm
and non detection may be tuned with the coefficient $\lambda^v$.

\begin{figure}
[h]
\begin{center}
\includegraphics[
height=1.9865in,
width=2.6385in
]%
{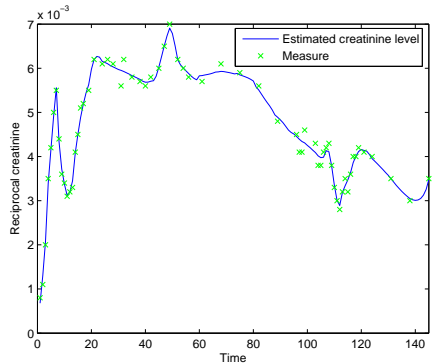}%
\caption{Measured (cross) and estimated (solid line) creatinine
level with a Rao-Blackwellized particle filter with 1000
particles.}
\label{ex2rbpfcreatinine}%
\end{center}
\end{figure}
\begin{figure}
[hh]
\begin{center}
\includegraphics[
height=2.0193in,
width=2.6826in
]%
{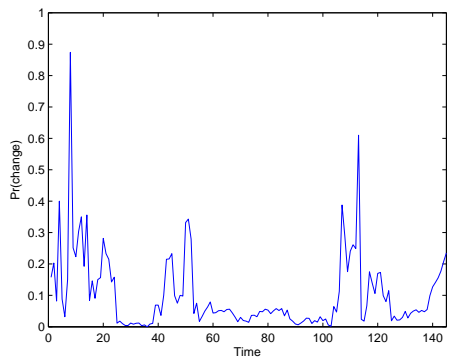}%
\caption{Posterior probability of a jump in the creatinine level
with the Rao-Blackwellized particle filter with 1000 particles.
For a threshold set to $0.5$, the creatinine level jumps are
detected at about times $8$ and $110$.}
\label{ex2rbpfsauts}%
\end{center}
\end{figure}

\section{Discussion}

\label{sec:discussion}In this section, we discuss several features of the
approach proposed.

\subsection{About Dirichlet Process-based modeling}

\label{sec:about-dirichl-proc} DPMs have several main advantages. Firstly,
sampling from the posterior distribution is made especially easy thanks to the
Polya urn scheme. Second, the discreteness of the distribution $\mathbb{G}$
enables straightforward estimation of the \textquotedblleft number of
components\textquotedblright, without requiring reversible jump-like
computational approaches. This discreteness has, however, some unexpected
effects on inferences, which are reported in~\cite{Petrone97}
and~\cite{Green01}. For example, the DP tends to favor a misbalance between
the size of the groups of latent variables associated to the same cluster, and
to concentrate the posterior distribution of the number of groups on a higher
value. Dirichlet Processes realize nevertheless an attractive trade-off
between versatile modeling properties and implementation advantages, which
explain their success in various contexts -- and our choice to use them in
this paper.





\subsection{About MCMC algorithms for DPMs}

As stated in~\cite{Neal00}, the ``single-site'' marginal algorithm used in this
paper may be stuck in a mode of the posterior: several noises samples
$\mathbf{v}_{t}$ (resp. $\mathbf{w}_{t}$) are associated to the same cluster
value $U_{j}^{v}$ for some $j$ in Eq.~(\ref{eq:stickbreak}) (resp.
$U_{j^{\prime}}^{w}$) -- in other words, there are many $t$'s such that
$\theta_{t}^{v}=U_{j}^{v}$ for some $j$ (resp $\theta_{t}^{w}=U_{j^{\prime}%
}^{w}$). Since the algorithm cannot change the value of $\theta_{t}^{v}$ for
more than one $\mathbf{v}_{t}$ simultaneously, changes to $\theta_{t}^{v}$
occur rarely, as they require passage through a low-probability intermediate
state in which noises $\mathbf{v}_{t}$ in the same group are not associated to
the same cluster. In alternative algorithms, such as those given
in~\cite{Neal00}, clusters are sampled in groups, which avoids this problem at
the expense of an increased computational cost. Nevertheless, we have
demonstrated empirically in Section~\ref{sec:applications} that our MCMC
scheme is indeed efficient in the applications presented.

\subsection{About the hyperparameter estimation in the MCMC algorithm}

As shown in the applications section, the estimation of the hyperparameter
$\alpha$ improves the overall state estimation. It also makes the convergence
of the Gibbs sampler faster. During the first iterations, the value of
$\alpha$ is high, and the sampler proposes new clusters more easily. This
enables efficient state space global exploration during the first iterations.
When the "good" clusters have been found, the value of $\alpha$ decreases, and
it eliminates useless clusters.

\subsection{About the convergence of the Rao-Blackwellized particle filter}

\label{sec:about-conv-rao} Because the DPMs $F^{v}$ and $F^{w}$ are static
(infinite-dimensional) parameters, the Rao-Blackwellized particle filter
suffers from an accumulation of errors over time. In other words, the particle
filter is not able to move cluster values $U_{j}^{v}$'s and $U_{j}^{w}$ after
they are initialized. This is a well known problem of \emph{static parameter}
estimation with particle filters. However, as the static component is not the
estimated cluster $\theta_{t}$ but its prior distribution $\mathbb{G}$, this
accumulation is less critical than with the estimation of true static parameters.

In Section~\ref{sec:SMC}, the hyperparameter vector $\phi$ is assumed fixed,
also because this is a static parameter. It could actually be estimated by
implementing one of the particle filtering approaches to static parameter
estimation. For example, the approaches
in~\cite{Liu01,Gilks01,Doucet03,Poyiadjis05} are based on either kernel
density methods, MCMC steps, or Maximum Likelihood. However, these algorithms
also have important drawbacks (error accumulation with time in $O(t^{2})$). An alternative solution consists of introducing
an artificial dynamic on the hyperparameters~\cite{Andrieu03} but it is not
applicable to our problem: we would then loose the Polya urn structure given
by Eq. (\ref{eq:Polya}).

\subsection{About related approaches}

\label{sec:about-relat-appr} Our model has some connections with Jump Linear
Systems (JLS)~\cite{Ackerson70,Akashi77}. In JLS, a discrete indicator
variable switches between a (known) fixed number of different (known) linear
Gaussian models with some (known) prior probability. Our model may be
interpreted as a JLS whose number of different models is unknown, mean vector
and covariance matrix of the linear Gaussian models are unknowns as well as
their prior probabilities. The model proposed in this paper can also be
generalized in the following manner. Denote $\underline{\theta}_{t}$ $=$
$\{F_{t},C_{t},H_{t},G_{t},\mu_{t}^{v},\Sigma_{t}^{v},\mu_{t}^{w},\Sigma
_{t}^{w}\}$ $=$ $\{F_{t},C_{t},H_{t},G_{t},\theta_{t}\}$ and $\underline
{\mathbb{G}}_{0}$ a prior distribution on $\underline{\theta}_{t}$. The
following general hierarchical model
\begin{equation}%
\begin{tabular}
[c]{l}%
$\underline{\mathbb{G}}\mathbb{\ }\sim\ DP(\underline{\mathbb{G}}_{0}%
,\alpha),$\\
$\underline{\theta}_{t}|\underline{\mathbb{G}}\sim\underline{\mathbb{G}},$\\
$\mathbf{x}_{t}|\underline{\theta}_{t},\mathbf{x}_{t-1}\ \sim\ \mathcal{N}%
(F_{t}\mathbf{x}_{t-1}+C_{t}\mathbf{u}_{t}+G_{t}\mu_{t}^{v},G_{t}\Sigma
_{t}^{v}G_{t}^{\prime}),$\\
$\mathbf{z}_{t}|\underline{\theta}_{t},\mathbf{x}_{t}\ \sim\ \mathcal{N}%
(H_{t}\mathbf{x}_{t}+\mu_{t}^{w},\Sigma_{t}^{w})$%
\end{tabular}
\ \ \ \ \ \ \
\end{equation}
has more flexibility than common JLS: the number of different switching models
is estimated, as well as the parameters of these models and their prior probabilities.

\subsection{About observability}

\label{sec:about-observability} In order for the observation noise
$\mathbf{w}_{t}$ pdf to be correctly estimated, some observability constraints
must be ensured. Indeed, the pair $(\widetilde{F},\widetilde{H})$ has to be
fully observable, that is, the observability matrix%
\begin{equation}
\left(
\begin{array}
[c]{c}%
\widetilde{H}\\
\widetilde{H}\widetilde{F}\\
\ldots\\
\widetilde{H}\widetilde{F}^{n_{x}+n_{z}-1}%
\end{array}
\right)
\end{equation}
must have rank $n_{x}+n_{z}$ (full rank), where $\widetilde{F}=\left(
\begin{array}
[c]{cc}%
F & 0_{n_{x}\times n_{z}}\\
0_{n_{z}\times n_{x}} & I_{n_{z}}%
\end{array}
\right)  $, $\widetilde{H}=\left(
\begin{array}
[c]{cc}%
H & I_{n_{z}}%
\end{array}
\right)  $, $n_{x}$ and $n_{z}$ are resp. the length of the state and
observation vectors.

\section{Conclusion}

In this paper, we have presented a Bayesian nonparametric model that enables
state and observation noise pdfs estimation, in a linear dynamic model. The
Dirichlet process mixture considered here is flexible and we have presented
two simulation-based algorithms based on Rao-Blackwellization which allows us
to perform efficiently inference. The approach has proven efficient in
applications -- in particular, we have shown that state estimation is possible
even though the dynamic and observation noises are of unknown pdfs. We are
currently investigating the following extensions of our methodology. First, it
would be of interest to consider nonlinear dynamic models. Second, it would be
important to develop time-varying Dirichlet process mixture models in cases
where the noise statistics are assumed to evolve over time.

\appendix

\subsection{Notations}

$\mu$ and $\Sigma$ are sampled from a Normal inverse Wishart distribution
$\mathbb{G}_{0}$ of hyperparameters $\mu_{0}$, $\kappa_{0}$, $\nu_{0}$,
$\Lambda_{0}$ if
\begin{align*}
\mu|\Sigma &  \sim\mathcal{N}(\mu_{0},\frac{\Sigma}{\kappa_{0}})\\
\Sigma^{-1}  &  \sim W(\nu_{0},\Lambda_{0}^{-1})
\end{align*}
where $W(\nu_{0},\Lambda_{0}^{-1})$ is the standard Wishart distribution.

\subsection{Backward forward recursion}

The quantities $P_{t-1|t}^{\prime-1}(\theta_{t:T})$ and $P_{t-1|t}^{\prime
-1}(\theta_{t:T})m_{t-1|t}^{\prime}(\theta_{t:T})$ defined in Section
\ref{sec:fixed-known-hyperp}\ always satisfy the following backward
information filter recursion.

\begin{enumerate}
\item Initialization\newline$P_{T|T}^{\prime-1}(\theta_{T})=H_{T}^{\mathsf{T}%
}(\Sigma_{T}^{w})^{-1}H_{T}$\newline$P_{T|T}^{\prime-1}(\theta_{T}%
)m_{T|T}^{\prime}(\theta_{T})=H_{T}^{\mathsf{T}}(\Sigma_{T}^{w})^{-1}%
(\mathbf{z}_{T}-\mu_{T}^{w})$

\item Backward recursion. For $t=T-1..1,$%
\begin{equation}
\Delta_{t+1}=\left[  I_{n_{v}}+B^{\mathsf{T}}(\theta_{t+1})P_{t+1|t+1}%
^{\prime-1}(\theta_{t+1|T})B(\theta_{t+1})\right]  ^{-1}%
\end{equation}%
\begin{equation}
P_{t|t+1}^{\prime-1}(\theta_{t+1:T})=F_{t+1}^{\mathsf{T}}P_{t+1|t+1}%
^{\prime-1}(\theta_{t+1:T})(I_{n_{x}}-B(\theta_{t+1})\Delta_{t+1}%
(\theta_{t+1:T})B^{\mathsf{T}}(\theta_{t+1})P_{t+1|t+1}^{\prime-1}%
(\theta_{t+1:T}))F_{t+1}\nonumber
\end{equation}%
\begin{align}
P_{t|t+1}^{\prime-1}(\theta_{t+1:T})m_{t|t+1}^{\prime}(\theta_{t+1:t})  &
=F_{t+1}^{\mathsf{T}}(\theta_{t+1})\times(I_{n_{x}}-P_{t+1|t+1}^{\prime
-1}(\theta_{t+1:T})B(\theta_{t+1})\Delta_{t+1}(\theta_{t+1:T})B^{\mathsf{T}%
}(\theta_{t+1}))\nonumber\\
&  \times P_{t+1|t+1}^{\prime-1}(\theta_{t+1:T})\left(  m_{t+1|t+1}^{\prime
}(\theta_{t+1:T})-\mathbf{u}_{t+1}^{\prime}(\theta_{t+1})\right)
\end{align}%
\begin{equation}
P_{t|t}^{\prime-1}(\theta_{t:T})=P_{t|t+1}^{\prime-1}(\theta_{t+1:T}%
)+H_{t}^{\mathsf{T}}(\Sigma_{t}^{w})^{-1}H_{t}%
\end{equation}%
\begin{equation}
P_{t|t}^{\prime-1}(\theta_{t:T})m_{t|t}^{\prime}(\theta_{t:T})=P_{t|t+1}%
^{\prime-1}(\theta_{t+1:T})m_{t|t+1}^{\prime}(\theta_{t+1:T})+H_{t}%
^{\mathsf{T}}(\Sigma_{t}^{w})^{-1}(\mathbf{z}_{t}-\mu_{t}^{w})
\end{equation}
where $B(\theta_{t})=G_{t}\times$chol$(\Sigma_{t}^{v})^{\mathsf{T}}$.
\end{enumerate}

For the Metropolis Hasting ratio, we need to compute the acceptance
probability only with a probability constant
\begin{equation}
p(\mathbf{z}_{1:T}|\theta_{1:T})\propto p(\mathbf{z}_{t}|\theta_{1:t}%
,\mathbf{z}_{1:t-1})\int_{\mathcal{X}}p(\mathbf{z}_{t+1:T}|\mathbf{x}%
_{t},\theta_{t+1:T})p(\mathbf{x}_{t}|\mathbf{z}_{1:t},\theta_{1:t}%
)d\mathbf{x}_{t}%
\end{equation}

If $\Sigma_{t|t}(\theta_{1:t})\neq0$ then it exists $\Pi_{t|t}(\theta_{1:t})$
and $Q_{t|t}(\theta_{1:t})$ such that $\Sigma_{t|t}(\theta_{1:t}%
)=Q_{t|t}(\theta_{1:t})\Pi_{t|t}(\theta_{1:t})Q_{t|t}^{T}(\theta_{1:t})$. The
matrices $Q_{t|t}(\theta_{1:t})$ and $\Pi_{t|t}(\theta_{1:t})$ are
straightforwardly obtained using the singular value decomposition of
$\Sigma_{t|t}(\theta_{1:t})$. Matrix $\Pi_{t|t}(\theta_{1:t})$ is a
$n_{t}\times n_{t},1\leq n_{t}\leq n_{x}$ diagonal matrix with the nonzero
eigenvalues of $\Sigma_{t|t}(\theta_{1:t})$ as elements. Then one has%
\begin{equation}%
\begin{tabular}
[c]{l}%
$p(\mathbf{z}_{1:T}|\theta_{1:T})\propto\mathcal{N}(\widehat{\mathbf{z}%
}_{t|t-1}(\theta_{1:t}),S_{t|t-1}(\theta_{1:t}))\left\vert \Pi_{t|t}%
(\theta_{1:t})Q_{t|t}^{\mathsf{T}}(\theta_{1:t})P_{t|t+1}^{\prime-1}%
(\theta_{t+1:T})Q_{t|t}(\theta_{1:t})+I_{n_{t}}\right\vert ^{-\frac{1}{2}}$\\
$\times\exp(-\frac{1}{2}\widehat{\mathbf{x}}_{t|t}^{\mathsf{T}}(\theta
_{1:t})P_{t|t+1}^{\prime-1}(\theta_{t+1:T})\widehat{\mathbf{x}}_{t|t}%
(\theta_{1:t})-2\widehat{\mathbf{x}}_{t|t}^{\mathsf{T}}(\theta_{1:t}%
)P_{t|t+1}^{\prime-1}(\theta_{t+1:T})m_{t|t+1}^{\prime}(\theta_{t+1:T})$\\
$-(m_{t|t+1}^{\prime}(\theta_{t+1:T})-\widehat{\mathbf{x}}_{t|t}(\theta
_{1:t}))^{\mathsf{T}}\times P_{t|t+1}^{\prime-1}(\theta_{t+1:T})A_{t|t}%
(\theta_{1:t})\times P_{t|t+1}^{\prime-1}(\theta_{t+1:T})(m_{t|t+1}^{\prime
}(\theta_{t+1:T})-\widehat{\mathbf{x}}_{t|t}(\theta_{1:t})))$%
\end{tabular}
\end{equation}
where
\begin{equation}
A_{t|t}(\theta_{1:t})=Q_{t|t}(\theta_{1:t})\left[  \Pi_{t|t}^{-1}(\theta
_{1:t})+Q_{t|t}^{\mathsf{T}}(\theta_{1:t})P_{t|t+1}^{\prime-1}(\theta
_{t+1:T})Q_{t|t}(\theta_{1:t})\right]  ^{-1}Q_{t|t}^{\mathsf{{T}}}%
(\theta_{1:t})
\end{equation}
The quantities $\widehat{\mathbf{x}}_{t|t}(\theta_{1:t})$, $\Sigma
_{t|t}(\theta_{1:t})$, $\widehat{\mathbf{z}}_{t|t-1}(\theta_{1:t})$ and
$S_{t|t-1}(\theta_{1:t})$ are, resp., the one-step ahead filtered estimate and
covariance matrix of $\mathbf{x}_{t}$, the innovation at time $t$, and the
covariance of this innovation. These quantities are provided by the Kalman
filter, the system being linear Gaussian conditional upon $\theta_{1:t}$.

\section*{Acknowledgment}

This work is partially supported by the Centre National de la Recherche
Scientifique (CNRS) and the Région Nord-Pas de Calais.

\bibliographystyle{IEEEtran}
\bibliography{bibnonparametricbayes}

\end{document}